\documentclass[12pt]{article}

\title{Improved estimation in a non-Gaussian parametric regression}

\author {Pchelintsev Evgeny\thanks{
Department of Mathematics and Mechanics, Tomsk State University,
Lenin str. 36, 634050 Tomsk, Russia, and Laboratoire de
Math\'ematiques Rapha\"{e}l Salem, UMR 6085 CNRS, Universit\'e de
Rouen, Avenue de l'Universit\'e BP.12, 76800 Saint Etienne du
Rouvray Cedex, France, e-mail: evgen-pch@yandex.ru}}

\usepackage{amssymb}
\usepackage{amsfonts}
\usepackage{amsmath}
\usepackage{amsthm}
\usepackage{enumerate}
\usepackage{multicol}
\usepackage{graphicx}

\newtheorem{theorem}{Theorem}[section]
\newtheorem{proposition}[theorem]{Proposition}
\newtheorem{lemma}[theorem]{Lemma}

\newtheorem{remark}{Remark}[section]
\newtheorem{corollary}[theorem]{Corollary}

\newcommand\cG{{\cal G}}

\newcommand\cL{{\cal L}}
\newcommand\cB{{\cal B}}

\newcommand\cN{{\cal N}}

\newcommand\cD{{\cal D}}

\def\bbr{{\mathbb R}}

\def\text#1{\hbox{#1}}
\def\proof{{\noindent \bf Proof. }}
\def\endproof{\hspace*{130mm}\mbox{\ $\qed$}}

\def\E{{\bf E}}

\def\Chi{{\bf 1}}
\def\d{\mathrm{d}}
\def\build #1_#2{\mathrel{\mathop{\kern 0pt #1}\limits_{#2}}}
\newcommand{\zs}[1]{{\mathchoice{#1}{#1}{\lower.25ex\hbox{$\scriptstyle#1$}}
{\lower0.25ex\hbox{$\scriptscriptstyle#1$}}}}

\begin{document}

\maketitle

\begin{abstract}
The paper considers the problem of estimating the parameters in a
continuous time regression model with a non-Gaussian noise of
pulse type. The noise is specified by the Ornstein--Uhlenbeck
process driven by the mixture of a Brownian motion and a compound
Poisson process. Improved estimates for the unknown regression
parameters, based on a special modification of the James--Stein
procedure with smaller quadratic risk than the usual least squares
estimates, are proposed. The developed estimation scheme is
applied for the improved parameter estimation in the discrete time
regression with the autoregressive noise depending on unknown
nuisance parameters.

\end{abstract}

\vspace*{5mm} \noindent {\bf Keywords} Non-Gaussian parametric
regression $\cdot$ Improved estimates $\cdot$ Pulse noise $\cdot$
Ornstein--Uhlenbeck process $\cdot$ Quadratic risk $\cdot$
Autoregressive noise

\vspace*{5mm} \noindent {\bf Mathematics Subject Classification
(2010)} 62H12 - 62M10

\bibliographystyle{plain}
\renewcommand{\columnseprule}{.1pt}

\section{Introduction}\label{sec:1}

In 1961, James and Stein proposed shrinkage estimates which
outperform in mean square accuracy the maximum likelihood
estimates in the problem of estimating the mean of a
multidimensional Gaussian vector with unit covariance matrix
(James and Stein 1961). This result stimulated the development of
the theory of improved estimation for different regression models
with dependent errors. Fourdrinier and Strawderman and Fourdrinier
and Wells solved the problem of improved parametric estimation in
the regression with dependent non-Gaussian observations under
spherically symmetric distributions of the noise (Fourdrinier and
Strawderman 1996; Fourdrinier and Wells 1994). Fourdrinier and
Pergamenshchikov and Konev and Pergamenchtchikov investigated the
problem of improved estimation in nonparametric setting
(Fourdrinier and Pergamenshchikov 2007; Konev and
Pergamenchtchikov 2010).

In this paper, we consider the problem of improved parametric
estimation for a continuous time regression with dependent
non-Gaussian noise of pulse type. The noise is specified by the
Ornstein--Uhlenbeck process which is known to capture important
distributional deviation from Gaussianity and to be appropriate
for modelling various dependence structures (Barndorff-Nielsen and
Shephard 2001).

Consider a regression model satisfying the equation
\begin{equation}\label{sec:1.1}
d y_{t}=\sum_{j=1}^p\theta_j\phi_j(t)d t+d \xi_{t},\quad 0\leq
t\leq n,
 \end{equation}
where $\theta=(\theta_1,...,\theta_p)'$\ (the notation $'$ holds
for transposition) is the vector of unknown parameters from a
compact set $\Theta\subset\mathbb{R}^{p}$;\ $(\phi_j(t))_{1\leq j
\leq p}$\ are one-periodic $[0,+\infty)\to\bbr$ functions,
orthonormal in the space $\cL_2[0,1]$.\ The noise
$(\xi_{t})_{t\geq 0}$\ in \eqref{sec:1.1} is assumed to be a
non-Gaussian Ornstein--Uhlenbeck process obeying the following
stochastic differential equation
\begin{equation}\label{sec:1.2}
d \xi_{t}=a \xi_{t}d t+d u_{t},
 \end{equation}
where $a\leq 0$,\ $(u_{t})_{t \geq 0}$\ is a Levy process which is
the mixture
\begin{equation}\label{sec:1.3}
u_{t}=\varrho_{1}w_{t}+\varrho_{2}z_{t}
 \end{equation}
of a standard Brownian motion $(w_t)_{t\geq 0}$\ and a compound
Poisson process $(z_t)_{t\geq 0}$\ defined as
\begin{equation}\label{sec:1.4}
z_{t}=\sum_{j=1}^{N_t} Y_{j},
 \end{equation}
where $(N_t)_{t\geq 0}$\ is a Poisson process with the intensity
$\lambda>0$,\ and $(Y_j)_{j\geq1}$\ is a sequence of i.i.d.
Gaussian random variables with parameters (0,1). The noise
parameters $a, \varrho_{1}, \varrho_{2}$\ and $\lambda$\ are
unknown.

The problem is to construct improved estimates for the unknown
vector parameter $\theta$\ on the basis of observations
$(y_t)_{0\leq t\leq n}$,\ which have higher precision as compared
with the least squares estimates (LSE).

It will be observed that the regression model \eqref{sec:1.1} is
conditionally Gaussian given the $\sigma$-algebra\
$\cG=\sigma\{N_{t}, t\geq 0\}$\ generated by the Poisson process.
In Section \ref{sec:3} it is shown that the problem of estimating
parameter $\theta$\ in \eqref{sec:1.1} can be reduced to that of
estimating the mean in a conditionally Gaussian distribution with
a random covariance matrix depending on the unknown nuisance
parameters. This enables one to construct a shrinkage estimates
for unknown parameters $(\theta_1,\ldots,\theta_p)$\ in
\eqref{sec:1.1}. The main result is given in Theorem
\ref{Th.sec:3.1} which claims that the proposed estimate has less
risk than the LSE.

In Section \ref{sec:2} we propose a special modification of
James--Stein procedure for solving the problem of estimating the
mean in a conditionally Gaussian distribution. This procedure
allows one to control the mean square accuracy of estimates. It is
shown (Theorem \ref{Th.sec:2.1}) that this estimate has less mean
square risk than the usual LSE.

The rest of the paper is organized as follows. In Section
\ref{sec:4} we apply Theorem \ref{Th.sec:2.1} to the problem of
parameter estimation in a discrete time regression under a
Gaussian autoregressive noise with unknown parameters. Appendix
contains some technical results.


\section{On improved estimation in a conditionally Gaussian regression}\label{sec:2}

In Section \ref{sec:3}, we will shown that the initial problem of
estimating the parameters $(\theta_1,\ldots,\theta_p)$\ in model
\eqref{sec:1.1} reduces to the following one. Suppose that the
observation $Y$\ is a $p$-dimensional random vector which obeys
the equation
\begin{equation}\label{sec:2.1}
Y=\theta+\xi,
\end{equation}
where $\theta$\ is an unknown constant vector parameter from some
compact set $\Theta\subset\mathbb{R}^{p}$,\ $\xi$\ is a
conditionally Gaussian random vector with zero mean and the
covariance matrix $\cD(\cG)$,\ i.e.
$Law(\xi|\cG)=\cN_{p}(0,\cD(\cG))$,\ where $\cG$\ is some fixed
$\sigma$-algebra.\

The problem is to estimate $\theta$.\

Consider a shrinkage estimate for $\theta$\ of the form
\begin{equation}\label{sec:2.2}
\theta^{*}=\left(1-\frac{c}{\|Y\|}\right)Y,
\end{equation}
where $c$\ is a positive constant which will be specified later.

The choice of this estimate \eqref{sec:2.2} is motivated by the
need to control the quadratic risk
\begin{equation*}
R(\theta,\tilde{\theta})=\E_{\theta}\|\theta-\tilde{\theta}\|^2
\end{equation*}
in the case of conditionally Gaussian model \eqref{sec:2.1}.

It will be observed that such control can not be provided by the
ordinary James--Stein estimate (obtained from \eqref{sec:2.2} by
the change of $\|Y\|$ to $\|Y\|^2$).

In order to get an explicit upper bound for the quadratic risk of
estimate \eqref{sec:2.2} we impose some conditions on the random
covariance matrix $\cD(\cG)$.\

Assume that

$(\bf{C_1})$\ {\sl There exists a positive constant $\lambda_*$,\
such that the minimal eigenvalue of matrix $\cD(\cG)$\ satisfies
the inequality}
\begin{equation*}
\lambda_{min}(\cD(\cG))\geq \lambda_*  \quad\mbox{a.s.}
\end{equation*}

$(\bf{C_2})$\ {\sl The maximal eigenvalue of the matrix
$\cD(\cG)$\ is bounded on some compact set
$\Theta\subset\mathbb{R}^{p}$\ from above in the sense that}
\begin{equation*}
\sup_{\theta\in\Theta}\E_\theta\lambda_{max}(\cD(\cG))\leq a^*,
\end{equation*}
{\sl where $a^*$\ is a known positive constant.}

Further we will introduce some notation. Let denote the difference
of the risks of estimate \eqref{sec:2.2} and LSE $\hat{\theta}=Y$\
as
\begin{equation*}
\Delta(\theta):=R(\theta^{*},\theta)-R(\hat{\theta},\theta).
\end{equation*}
We will need also the following constant
\begin{equation*}
\gamma_p=\dfrac{\sum_{j=0}^{p-2}2^{\frac{j-1}{2}}(-1)^{p-j}\mu^{p-1-j}
\Gamma\left(\frac{j+1}{2}\right)-(-\mu)^pI(\mu)}{2^{p/2-1}\Gamma\left(\frac{p}{2}\right)d},
\end{equation*}
where $\mu=d/\sqrt{a^*}$,\
\begin{equation*}
I(a)=\int_0^\infty\frac{\exp(-r^2/2)}{a+r}dr\, \quad\mbox{and}
\quad d=\sup\{\|\theta\|: \theta\in\Theta\}.
\end{equation*}

\begin{theorem}\label{Th.sec:2.1}

Let the noise $\xi$\ in \eqref{sec:2.1} have a conditionally
Gaussian distribution $\cN_{p}(0,\cD(\cG))$\ and its covariance
matrix $\cD(\cG)$\ satisfy conditions $(\bf{C_1}), (\bf{C_2})$\
with some compact set $\Theta\subset\mathbb{R}^{p}$.\ Then the
estimator \eqref{sec:2.2} with $c=(p-1)\lambda_{*}\gamma_p$\
dominates the LSE $\hat{\theta}_{ML}$\ for any $p\geq 2$,\ i.e.
\begin{equation*}
\sup_{\theta\in\Theta}\Delta(\theta)\leq-[(p-1)\lambda_{*}\gamma_p]^2.
\end{equation*}

\end{theorem}

\proof First we will find the lower bound for the random variable
$\|Y\|^{-1}$.\

\begin{lemma}\label{Le.sec:2.2}
Under the conditions of Theorem 2.1
\begin{equation*}
\inf_{\theta\in\Theta}\E_{\theta}\frac{1}{\|Y\|}\geq\gamma_p.
\end{equation*}
\end{lemma}
The proof of lemma is given in the Appendix.

In order to obtain the upper bound for $\Delta(\theta)$\ we will
adjust the argument in the proof of Stein's lemma (James and Stein
1961) to the model \eqref{sec:2.1} with a random covariance
matrix.

We represent the risks of LSE and of \eqref{sec:2.2} as
\begin{gather*}
R(\hat{\theta},\theta)=\E_{\theta}\|\hat{\theta}-\theta\|^{2}
=\E_{\theta}(\E\|\hat{\theta}-\theta\|^{2}|\cG)
=\E_{\theta}tr\cD(\cG);\\[2mm]
R(\theta^{*},\theta)=R(\hat{\theta},\theta)+\E_{\theta}[\E((g(Y)-1)^{2}\|Y\|^{2}|\cG)]\\[2mm]
+2\sum_{j=1}^{p}\E_{\theta}[\E((g(Y)-1)Y_{j}(Y_{j}-\theta_j)|\cG)],
\end{gather*}
where $g(Y)=1-c/\|Y\|$.

Denoting $f(Y)=(g(Y)-1)Y_{j}$\ and applying the conditional
density of distribution of a vector $Y$\ with respect to
$\sigma$-algebra $\cG$
\begin{equation*}
p_{Y}(x|\cG)=\frac{1}{(2\pi)^{p/2}\sqrt{\det\cD(\cG)}}
\exp\left(-\frac{(x-\theta)'\mathcal{D}^{-1}(\cG)(x-\theta)}{2}\right),
 \end{equation*}
one gets
\begin{equation*}
I_{j}:=\E(f(Y)(Y_{j}-\theta_j)|\cG)=\int_{\mathbb{R}^p}f(x)(x-\theta_j)p_{Y}(x|\cG)
d x, \quad j=\overline{1,p}.
\end{equation*}

Making the change of variable
$u=\mathcal{D}^{-1/2}(\cG)(x-\theta)$\ and assuming
$\tilde{f}(u)=f(\mathcal{D}^{1/2}(\cG)u+\theta)$,\ one finds that
\begin{equation*}
I_{j}=\frac{1}{(2\pi)^{p/2}}\sum_{l=1}^{p}\langle
\mathcal{D}^{1/2}(\cG)\rangle_{jl}\int_{\mathbb{R}^{p}}\tilde{f}(u)u_{l}\exp\left(-\frac{\|u\|^{2}}
{2}\right) d u, \quad j=\overline{1,p},
 \end{equation*}
where $\langle A\rangle_{ij}$\ denotes the $(i,j)$-th\ element of matrix $A$.\ These quantities can be written as
\begin{equation*}
I_{j}=\sum_{l=1}^{p}\sum_{k=1}^{p}\E(<\mathcal{D}^{1/2}(\cG)>_{jl}
<\mathcal{D}^{1/2}(\cG)>_{kl}\frac{\partial f}{\partial
u_k}(u)|_{u=Y}|\cG), \quad j=\overline{1,p}.
\end{equation*}

Thus, the risk for an estimator \eqref{sec:2.2} takes the form
\begin{gather*}
R(\theta^{*},\theta)=R(\hat{\theta},\theta)+\E_{\theta}((g(Y)-1)^{2}\|Y\|^{2})\\
+2\E_{\theta}\left(\sum_{j=1}^{p}\sum_{l=1}^{p}\sum_{k=1}^{p}<\mathcal{D}^{1/2}(\cG)>_{jl}
<\mathcal{D}^{1/2}(\cG)>_{kl}\frac{\partial}{\partial
u_{k}}[(g(u)-1)u_j]|_{u=Y}\right).
\end{gather*}

Therefore, one has
\begin{equation*}
R(\theta^{*},\theta)=R(\hat{\theta},\theta)+\E_{\theta}W(Y),
\end{equation*}
where
\begin{equation*}
W(z)=c^{2}+2c\frac{z'\cD(\cG)z}{\|z\|^3}-2tr\mathcal{D}(\textit{G})c\frac{1}{\|z\|}.
\end{equation*}
This implies that
\begin{equation*}
\Delta(\theta)=\E_{\theta}W(Y).
\end{equation*}

Since $z'Az\leq\lambda_{max}(A)\|z\|^2$,\ one comes to the inequality
\begin{equation*}
\Delta(\theta)\leq
c^2-2c\E_{\theta}\frac{tr\cD(\cG)-\lambda_{max}(\cD(\cG))}{\|Y\|}.
\end{equation*}
From here, it follows that
\begin{equation*}
\Delta(\theta)\leq
c^2-2c\sum_{i=2}^p\E_{\theta}\frac{\lambda_{i}(\cD(\cG))}{\|Y\|}.
\end{equation*}
Taking into account the condition $(\bf{C_1})$\ and the Lemma
~\ref{Le.sec:2.2}, one obtains
\begin{equation*}
\Delta(\theta)\leq c^2-2(p-1)\lambda_{*}\gamma_p c=:\phi(c).
\end{equation*}
Minimizing the function $\phi(c)$\ with respect to $c$,\ we come
to the desired result, i.e.
\begin{equation*}
\Delta(\theta)\leq -[(p-1)\lambda_{*}\gamma_p]^2.
\end{equation*}

Hence Theorem~\ref{Th.sec:2.1}.

\endproof

\begin{corollary}\label{Co.sec:2.3}
Let  in \eqref{sec:2.1} the noise $\xi\sim\cN_{p}(0,D)$\ with the
positive definite non random covariance matrix $D>0$\ and
$\lambda_{min}(D)\geq\lambda_*>0$.\ Then the estimator
\eqref{sec:2.2} with $c=(p-1)\lambda_*\gamma_p$\ dominates the LSE
for any $p\geq 2$\ and compact set $\Theta\subset\mathbb{R}^{p}$,\
i.e.
\begin{equation*}
\sup_{\theta\in\Theta}\Delta(\theta)\leq-[(p-1)\lambda_{*}\gamma_p]^2.
\end{equation*}

\end{corollary}

\begin{remark}
Note that if $D=\sigma^2I_p$\ then
\begin{equation*}
\sup_{\theta\in\Theta}\Delta(\theta)\leq-[(p-1)\sigma^2\gamma_p]^2.
\end{equation*}
\end{remark}

\begin{corollary}\label{Co.sec:2.4}
If $\xi\sim\cN_{p}(0,I_p)$ and $\theta=0$\ in model
\eqref{sec:2.1} then the risk of estimate \eqref{sec:2.2} is given
by the formula
\begin{equation}\label{sec:1.6}
R(0,\theta^*)=p-\left[\frac{(p-1)\Gamma((p-1)/2)}{\sqrt{2}\Gamma(p/2)}\right]^2=:r_p.
\end{equation}
\end{corollary}

By applying the Stirling's formula for the Gamma function
\begin{equation*}
\Gamma(x)=\sqrt{2\pi}x^{x-1/2}\exp(-x)\left(1+o(1)\right)
\end{equation*}
one can check that $r_p\rightarrow 0.5$\ as $p\rightarrow\infty$.\
The behavior of the risk \eqref{sec:1.6} for small values of $p$\
is shown in Fig.1. It will be observed that in this case the risk
of the James--Stein estimate $\hat{\theta}_{JS}$\ remains constant
for all $p\geq3$,\ i.e.
\begin{equation*}
R(0,\hat{\theta}_{JS})=2
\end{equation*}
and the risk of the LSE $\hat{\theta}$\ is equal to $p$\ and tends
to infinity as $p\rightarrow\infty$.\

\begin{figure}[t]
\centering
\includegraphics[width=0.9\textwidth]{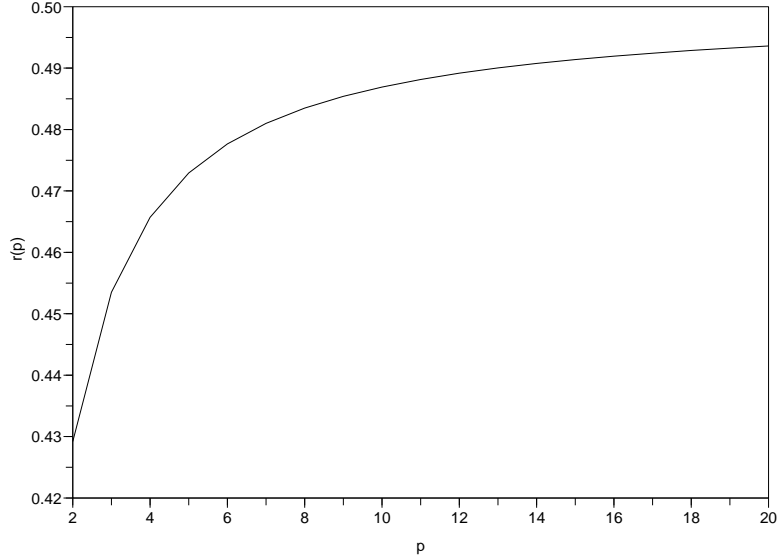}
\caption{Risk of $\theta^*$\ at $\theta=0$.\ }
\label{Risk_figure}
\end{figure}

\section{Improved estimation in a non-Gaussian Ornstein--Uhlenbeck--Levy regression model}\label{sec:3}

In this section we use the estimate \eqref{sec:2.2} to a
non-Gaussian continuous time regression model to construct an
improved estimate of the unknown vector parameter $\theta$.\ To
this end we reduce first the initial continuous time regression
model \eqref{sec:1.1} to a discrete time model of the form
\eqref{sec:2.1} with a conditionally Gaussian noise.

A commonly used estimator of an unknown vector $\theta$\ in model
\eqref{sec:1.1} on the basis of observations $(y_t)_{0\leq t\leq
n}$\ is the LSE
$\hat{\theta}=(\hat{\theta}_1,\ldots,\hat{\theta}_p)'$\ with the
components
\begin{equation*}
\hat{\theta_j}=\frac{1}{n}\int_0^n\phi_j(t)dy_t,\quad j=\overline{1,p}.
\end{equation*}
From here and \eqref{sec:1.1}, one has
\begin{equation}\label{sec:3.1}
\hat{\theta}=\theta+n^{-1/2}\zeta(n),
\end{equation}
where $\zeta(n)$\ is the random vector with the coordinates
\begin{equation*}
\zeta_j(n)=n^{-1/2}\int_0^n\phi_j(t)d\xi_t.
\end{equation*}
Note that the vector $\zeta(n)$\ has a conditionally Gaussian
distribution with a zero mean and conditional covariance matrix
$V_{n}(\cG)=cov(\zeta(n),\zeta(n)'|\cG)$\ with the elements
\begin{equation*}
v_{ij}(n)=\E(\zeta_{i}(n)\zeta_{j}(n)\mid\cG).
\end{equation*}

Thus the initial problem of estimating parameter $\theta$\ in
\eqref{sec:1.1} can be reduced to the that of estimating parameter
$\theta$\ in a conditionally Gaussian regression model
\eqref{sec:3.1}.

\begin{theorem}\label{Th.sec:3.1}
Let the regression model be given by the equations
\eqref{sec:1.1}--\eqref{sec:1.4}, $\varrho_{1}>0$.\ Then, for any
$n\geq 1$\ and $p\geq 2$,\ the estimator of $\theta$\
\begin{equation*}
\theta^*=\left(1-\frac{\varrho_{1}^{2}(p-1)\gamma_p}{n\Vert\hat{\theta}\Vert}\right)\hat{\theta},
\end{equation*}
dominates the LSE $\hat{\theta}$:\
\begin{equation*}
\sup_{\theta\in\Theta}\Delta(\theta)\leq -\left[\frac{\varrho_{1}^{2}(p-1)\gamma_p}{n}\right]^2.
\end{equation*}

\end{theorem}

To prove this theorem one can apply Theorem \ref{Th.sec:2.1}. To
this end it suffices to check conditions $(\bf{C_1})$,\
$(\bf{C_2})$\ on the matrix $V_{n}(\cG)$.\ The verification of
conditions $(\bf{C_1})$\ and $(\bf{C_2})$\ is given in the
Appendix.

\section{Improved estimation in an autoregression}\label{sec:4}

In this section we consider the problem of improved estimating the
unknown mean of a multivariate normal distribution when the
dispersion matrix is unknown and depends on some nuisance
parameters. The models of autoregressive type are widely used in
time series analysis (Anderson 1994; Brockwell and Davis 1991).

Let the noise $\xi=(\xi_1,\ldots,\xi_p)'$\ in \eqref{sec:2.1}  be
described by a Gaussian autoregression process
\begin{equation}\label{sec:4.1}
\xi_k=a\xi_{k-1}+\varepsilon_k,\ k=\overline{1,p},
\end{equation}
where $|a|<1$,\ $\E\xi_0=0$\ and
$\varepsilon_1,\ldots,\varepsilon_p$\ are independent Gaussian
(0,1) random variables. Assume that the parameter $a$\ in
\eqref{sec:4.1} is unknown and belongs to interval
$[-\alpha,\alpha]$,\ where $0<\alpha<1$\ is known number.

It is easy to check that the covariance of the noise $\xi$\ has the form

\begin{equation*}
D(a)=\frac{1}{1-a^2} \left(
\begin{array}{llll}
1 & a &\ldots& a^{p-1} \\[2mm]
     a & 1 &\ldots& a^{p-2} \\[2mm]
     &  \ddots & &   \\[2mm]
a^{p-1} & a^{p-2} &\ldots& 1
\end{array}
\right)
\end{equation*}

\begin{proposition}\label{Pr.sec:4.1}
Let the noise $\xi$ in \eqref{sec:2.1} be specified by equation
\eqref{sec:4.1} with $a\in [-\alpha,\alpha]$.\ Then, for any
$p>1/(1-\alpha)^2$,\ the LSE is dominated by the estimate
\begin{equation*}
\theta^*=\left(1-\left(p-\frac{1}{(1-\alpha)^2}\right)\frac{\gamma_p}{\|Y\|}\right)Y
\end{equation*}
in the sense that
\begin{equation*}
\sup_{\theta\in\Theta}\Delta(\theta)\leq -\left(p-\frac{1}{(1-\alpha)^2}\right)^2\gamma_p^2.
\end{equation*}

\end{proposition}

\proof We note that $tr D(a)=p/(1-a^2)$.\ Now we will estimate of
the maximal eigenvalue of matrix $D(a)$.\ From the definition
\begin{equation*}
\lambda_{max}(D(a))=\sup_{\Vert z\Vert=1}z'D(a)z
\end{equation*}
one has
\begin{gather*}
z'D(a)z=\sum_{i=1}^{p}\sum_{j=1}^{p}<D(a)>_{ij}z_iz_j=\frac{1}{1-a^2}\left(1+2\sum_{i=1}^{p-1}\sum_{j=1}^{p-i}a^jz_iz_{j+i}\right)\\
=\frac{1}{1-a^2}\left(1+2\sum_{j=1}^{p-1}a^j\sum_{i=1}^{p-j}z_jz_{i+j}\right).
\end{gather*}
Applying the Cauchy--Bunyakovskii inequality yields
\begin{equation*}
\lambda_{max}(D(a))\leq\frac{1}{1-\alpha^2}\left(1+2\sum_{j=1}^{\infty}\alpha^j\right)=\frac{1}{(1-\alpha)^2}.
\end{equation*}
Thus,
\begin{equation*}
tr D(a)-\lambda_{max}(D(a))\geq p-\frac{1}{(1-\alpha)^2}.
\end{equation*}
By applying Theorem \ref{Th.sec:2.1} we come to the assertion of
Proposition \ref{Pr.sec:4.1}.

\endproof

\section{Appendix}\label{sec:5}

\subsection{Proof of the Lemma \ref{Le.sec:2.2}.}

\proof From \eqref{sec:2.1}, one has
\begin{equation*}
J=\E_{\theta}\frac{1}{\|Y\|}=\E_{\theta}\frac{1}{\|\theta+\xi\|}\geq\E_{\theta}\frac{1}{d+\|\xi\|}.
\end{equation*}
Taking the repeated conditional expectation and noting that the
random vector $\xi$\ is conditionally Gaussian with zero mean, one
gets
\begin{equation*}
J\geq\E_\theta\frac{1}{(2\pi)^{p/2}\sqrt{det\cD(\cG)}}
\int_{\mathbb{R}^p}\frac{\exp(-x'\cD(\cG)^{-1}x/2)}{d+\|x\|}dx.
\end{equation*}
Making the change of variable $u=\cD(\cG)^{-1/2}x$\ and applying
the estimation $u'\cD(\cG)u\leq\lambda_{max}(\cD(\cG))\|u\|^2$,\
we find
\begin{equation*}
J\geq\frac{1}{(2\pi)^{p/2}}
\int_{\mathbb{R}^p}\frac{\exp(-\|u\|^2/2)}{d+\sqrt{\lambda_{max}(\cD(\cG))}\|u\|}du.
\end{equation*}

Further making the spherical changes of the variables yields

\begin{equation*}
J\geq\frac{1}{2^{p/2-1}\Gamma(p/2)}\E_\theta\int_{0}^\infty\frac{r^{p-1}\exp(-r^2/2)}{d+\sqrt{\lambda_{max}(\cD(\cG))}r}dr.
\end{equation*}

From here by applying the Jensen and Cauchy--Bunyakovskii
inequalities and by the condition $(\mathbf{C_2})$,\ we obtain
\begin{equation*}
J\geq\frac{\mu}{2^{p/2-1}\Gamma(p/2)d}
\int_{0}^\infty\frac{r^{p-1}\exp(-r^2/2)}{\mu+r}dr=\gamma_p.
\end{equation*}

This leads to the assertion of Lemma \ref{Le.sec:2.2}.

\endproof

\subsection{The verification of the conditions $(\bf{C_1})$\ and
$(\bf{C_2})$\ on the matrix $V_{n}(\cG)$.\ }

Now we establish some properties of a stochastic integral
\begin{equation*}
I_\zs{t}(f)=\int^{t}_\zs{0}\, f_\zs{s} d \xi_\zs{s},\quad 0\le
t\le n
\end{equation*}
with respect to the process \eqref{sec:1.2}. We will need some
notations. Let us denote
\begin{equation*}
\varepsilon_\zs{f}(t)=a\int^{t}_\zs{0}\,\exp\{a(t-v)\}\,f(v)\,(1+\exp(2av))\,d
v\,,
\end{equation*}
where $f$ is $[0,+\infty)\to\bbr$ function integrated on any
finite interval. We introduce also the following transformation
\begin{equation*}
\tau_\zs{f,g}(t)= \frac{1}{2} \int^{t}_\zs{0} \left(2 f(s)g(s) +
\varepsilon^{*}_\zs{f,g}(s) \right)\, d s
\end{equation*}
of square integrable $[0,+\infty)\to\bbr$ functions $f$ and $g$.
Here
\begin{equation*}
\varepsilon^{*}_\zs{f,g}(t)=f(t)\varepsilon_\zs{g}(t) +
\varepsilon_\zs{f}(t)g(t)\,.
\end{equation*}

\begin{proposition}\label{Pr.sec:5.1}
If $f$ and $g$ are functions from $\cL_\zs{2}[0,n]$ then
\begin{equation}\label{sec:5.1}
\E\, I_\zs{t}(f)I_\zs{t}(g)=\varrho^{*}\,\tau_\zs{f,g}(t)
\end{equation}
where $\varrho^{*}=\varrho_{1}^{2}+\lambda\varrho_{2}^{2}$.\
\end{proposition}
\proof Noting that the process $I_\zs{t}(f)$ satisfies the
stochastic equation
$$
d I_\zs{t}(f)=a f(t)\xi_\zs{t}d t+f(t) d u_\zs{t}\,,\quad
I_\zs{0}(f)=0\,,
$$
and applying the Ito formula, one obtains \eqref{sec:5.1}. Hence
Proposition~\ref{Pr.sec:5.1}.

\endproof

\begin{corollary}
If $f$ is function from $\cL_\zs{2}[0,n]$ then
\begin{equation}\label{sec:5.2}
\E\, I_\zs{n}^2(f)\leq 3\varrho^{*}\,\int_0^nf^2(t)dt.
\end{equation}
\end{corollary}

Further, for an integrated $[0,+\infty)\to\bbr$ function $f$, we
define the function
$$
L_\zs{f}(x,z)=a \exp(a x)\,\left( f(z)+a\int^{x}_\zs{0}\,\exp(a
v)\,f(v+z)\,d v \right).\
$$
Let $(T_l)_{l\geq 1}$\ be the jump times of the Poisson process
$(N_t)_{t\geq 0}$,\ i.e.
\begin{equation*}
T_{l}=\inf\{t\geq 0: N_{t}=l\}.
\end{equation*}

\medskip

\begin{proposition}\label{Pr.sec:5.2}
Let $f$ and $g$  be  bounded left-continuous
$[0,\infty)\times\Omega\to\bbr$ functions
 measurable with respect to $\cB[0,+\infty)\bigotimes \cG$
(the product $\sigma$ algebra created by $\cB[0,+\infty)$ and
$\cG$). Then
$$
\E\left( I_\zs{t}(f) |\cG \right) =0
$$
and
\begin{multline*}
 \E\left( I_\zs{t}(f)\, I_\zs{t}(g) |\cG
\right)=\varrho^2_\zs{1}\tau_\zs{f,g}(t)+\varrho^2_\zs{2}\sum_{l\geq 1}f(T_l)g(T_l)\Chi_{(T_{l}\leq t)} \\
+\varrho_{2}^{2}\sum_{l\geq1}\int_{0}^{t}\left(f(s)L_{g}(s-T_{l},T_{l})+g(s)L_{f}(s-T_{l},T_{l})\right)\Chi_{(T_{l}\leq
s)}d s.
\end{multline*}

\end{proposition}
\proof By the Ito formula one has
\begin{align}\nonumber
I_\zs{t}(f)\, I_\zs{t}(g)&= \int^{t}_\zs{0}\,(\varrho^{2}_\zs{1}
f(s)g(s)
+a (f(s) I_\zs{s}(g) + g(s) I_\zs{s}(f))\xi_\zs{s})\d s\\[2mm] \nonumber
&+\varrho^{2}_\zs{2}\, \sum_\zs{l\ge
1}\,f(T_\zs{l})\,g(T_\zs{l})\,Y^{2}_\zs{l} \Chi_\zs{\{T_\zs{l}\le
t\}}
\\[2mm]
\nonumber& +\int^{t}_\zs{0}\, (f(s)I_\zs{s-}(g) +
g(s)I_\zs{s-}(f)))\d u_\zs{s}\,.
\end{align}
Taking the conditional expectation $\E\left(\cdot|\cG\right)$, on
the set $\{T_\zs{l}>t\}$, yields
\begin{align*}
\E \left(I_\zs{t}(f)\, I_\zs{t}(g)|\cG\right)&=
\int^{t}_\zs{0}\,\varrho^{2}_\zs{1} f(s)g(s)\d s
+\varrho^{2}_\zs{2}\, \sum_\zs{l\ge 1}\,f(T_\zs{l})\,g(T_\zs{l})\,
\,\Chi_\zs{\{T_\zs{l}\le t\}}
\\[2mm]
&+a \int^{t}_\zs{0}\, \left( f(s) \E(I_\zs{s}(g) \xi_\zs{s}|\cG) +
g(s) \E(I_\zs{s}(f)\xi_\zs{s}|\cG) \right) \d s\,.
\end{align*}

It is easy to check that
$$
a\E(I_\zs{t}(f)\xi_\zs{t}|\cG)=
\frac{\varrho^{2}_\zs{1}}{2}\varepsilon_\zs{f}(t)
+\varrho^{2}_\zs{2}\,\sum_\zs{j\ge
1}\,L_\zs{f}(t-T_\zs{j},T_\zs{j})\, \Chi_\zs{\{T_\zs{j}\le t\}}.
$$
From here one comes to the desired equality. Hence
Proposition~\ref{Pr.sec:5.2}.

\endproof

Thus, in view of $\zeta_j(n)=n^{-1/2}I_n(\phi_j)$ and
Proposition~\ref{Pr.sec:5.2} the elements of covariance matrix
$V_{n}(\cG)$\ can be written as
\begin{multline}\label{sec:5.3}
v_{ij}(n)=\frac{\varrho_{1}^{2}}{n}\int_{0}^{n}\phi_i(t)\phi_j(t)dt\\
+\frac{\varrho_{1}^{2}}{2n}\int_{0}^{n}\left(\phi_i(t)\varepsilon_{\phi_j}(t)+\phi_j(t)\varepsilon_{\phi_i}(t)\right)dt
+\frac{\varrho_{2}^{2}}{n}\sum_{l\geq 1}\phi_i(T_l)\phi_j(T_l)\Chi_{(T_{l}\leq n)}\\
+\frac{\varrho_{2}^{2}}{n}\sum_{l\geq1}\int_{0}^{n}\left(\phi_i(t)L_{\phi_j}(t-T_{l},T_{l})+\phi_j(t)L_{\phi_i}(t-T_{l},T_{l})\right)\Chi_{(T_{l}\leq
t)}dt.
\end{multline}

\begin{lemma}\label{Lem.sec:5.4}
Let $(\xi_{t})_{t\geq 0}$\ be defined by \eqref{sec:1.2} with
$a\leq 0$.\ Then a matrix $V_{n}(\cG)=(v_{ij}(\phi))_{1\leq
i,j\leq p}$\ with elements defined by \eqref{sec:5.3}, satisfy the
following inequality a.s.
\begin{equation*}
\inf_{n\geq 1}\inf_{\|z\|=1}z'V_{n}(\cG)z\geq\varrho_{1}^{2}.
\end{equation*}
\end{lemma}

\proof Notice that by \eqref{sec:5.3} one can the matrix
$V_{n}(\cG)$\ present as
\begin{equation*}
V_{n}(\cG)=\varrho_{1}^{2}I_p+F_n+B_{n}(\cG),
\end{equation*}
where $F_n$\ is non random matrix with elements
\begin{equation*}
f_{ij}(n)=\frac{\varrho_{1}^{2}}{2n}\int_{0}^{n}\left(\phi_i(t)\varepsilon_{\phi_j}(t)+\phi_j(t)\varepsilon_{\phi_i}(t)\right)dt
\end{equation*}
and $B_{n}(\cG)$\ is a random matrix with elements
\begin{multline*}
b_{ij}(n)=\frac{\varrho_{2}^{2}}{n}\sum_{l\geq
1}[\phi_i(T_l)\phi_j(T_l)\Chi_{(T_{l}\leq n)}\\
+\int_{0}^{n}\left(\phi_i(t)L_{\phi_j}(t-T_{l},T_{l})+\phi_j(t)L_{\phi_i}(t-T_{l},T_{l})\right)\Chi_{(T_{l}\leq
t)}d t ].
 \end{multline*}
This implies that
\begin{equation*}
z'V_{n}(\cG)z=\varrho_{1}^{2}z'z+z'F_nz+z'B_{n}(\cG)z\geq
\varrho_{1}^{2}z'z,
\end{equation*}
therefore
\begin{equation*}
\inf_{\|z\|=1}z'V_{n}(\cG)z\geq\varrho_{1}^{2}
\end{equation*}
and we come to the assertion of Lemma \ref{Lem.sec:5.4}.

\endproof

\begin{lemma}\label{Lem.sec:5.5}
Let $(\xi_{t})_{t\geq 0}$\ be defined by \eqref{sec:1.2} with
$a\leq 0$.\ Then a maximal eigenvalue of the matrix
$V_{n}(\cG)=(v_{ij}(n))_{1\leq i,j\leq p}$\ with elements defined
by \eqref{sec:5.3}, satisfy the following inequality
\begin{equation*}
\sup_{n\geq1}\sup_{\theta\in\Theta}\E_\theta\lambda_{max}(V_{n}(\cG))\leq
3 p\varrho^*.
\end{equation*}
where $\varrho^*=\varrho_{1}^{2}+\lambda\varrho_{2}^{2}$.\

\end{lemma}

\proof We note that
\begin{equation*}
\E_\theta\lambda_{max}(V_{n}(\cG))\leq\E_\theta
tr(V_{n}(\cG))=\sum_{j=1}^p\E_\theta
\zeta_j^2(n)=\frac{1}{n}\sum_{j=1}^p\E_\theta I_n^2(\phi_j).
\end{equation*}
Applying \eqref{sec:5.2} and $\int_0^n\phi_j^2(t)dt=n$, we obtain
the desired inequality.

Hence Lemma \ref{Lem.sec:5.5}.

\endproof

Thus the matrix $V_{n}(\cG)$\ is positive definite and satisfies
for any compact set $\Theta\subset\mathbb{R}^{p}$,\ the conditions
$(\bf{C_1})$\ and $(\bf{C_2})$\ with $\lambda_*=\varrho_1^2$\ and
$a^*=3p\varrho^*$.\


\begin{thebibliography}{}

\end{thebibliography}


\medskip

\begin{thebibliography}{100}

\bibitem{1}
Anderson T.W. (1994) The Statistical Analysis of Time Series.
Wiley, New York, London, Sydney, Toronto

\bibitem{2}
Barndorff-Nielsen O.E., Shephard N. (2001) Non-Gaussian
Ornstein-Uhlenbeck-based models and some of their uses in
financial mathematics. J. Royal Stat. Soc. B 63:167-241

\bibitem{3}
Brockwell P.J., Davis R.A. (1991) Time Series: Theory and Methods,
2nd edn. Springer, New York

\bibitem{4}
Fourdrinier D., Pergamenshchikov S. (2007) Improved selection
model method for the regression with dependent noise, Ann. of the
Inst. of Statist. Math. 59(3):435-464

\bibitem{5}
Fourdrinier D., Strawderman W.E. (1996) A paradox concerning
shrinkage estimators: should a known scale parameter be replaced
by an estimated value in the shrinkage factor? J. Multivariate
Anal. 59(2):109-140

\bibitem{6}
Fourdrinier D., Wells M.T. (1994) Comparaison de proc\'edures de
s\'election d'un mod\`{e}le de r\'egression: une approche
d\'ecisionnelle, Comptes Rendus de l'Acad\'emie des Sciences
Paris, v.319, s\'erie I:865-870

\bibitem{7}
James W., Stein C. (1961) Estimation with quadratic loss, in:
Proceedings of the Fourth Berkeley Symposium on Mathematics
Statistics and Probability, Vol. 1, University of California
Press, Berkeley, p.361-380

\bibitem{8}
Konev V., Pergamenchtchikov S. (2010) General model selection of a
periodic regression with a Gaussian noise. An. Inst. Stat. Math.
62:1083-1111


\end{thebibliography}
\end{document}